\documentclass{article}
\usepackage[utf8]{inputenc}
\usepackage{
 amsmath,
 amsxtra,
 amsthm,
 amssymb,
 etex,
 mathrsfs,
  }
\usepackage[all]{xy}

\newtheorem{theorem}{Theorem}[section]
\newtheorem{lemma}[theorem]{Lemma}
\newtheorem{conjecture}[theorem]{Conjecture}
\newtheorem{proposition}[theorem]{Proposition}
\newtheorem{corollary}[theorem]{Corollary}

\theoremstyle{definition}
\newtheorem{defn}[theorem]{Definition}
\newtheorem{remark}[theorem]{Remark}

\newcommand{\bd}{\begin{defn}}
\newcommand{\ed}{\end{defn}}
\newcommand{\bl}{\begin{lemma}}
\newcommand{\el}{\end{lemma}}
\newcommand{\bp}{\begin{proposition}}
\newcommand{\ep}{\end{proposition}}
\newcommand{\bt}{\begin{theorem}}
\newcommand{\et}{\end{theorem}}
\newcommand{\bc}{\begin{corollary}}
\newcommand{\ec}{\end{corollary}}
\newcommand{\br}{\begin{remark}}
\newcommand{\er}{\end{remark}}
\newcommand{\ba}{\begin{array}}
\newcommand{\ea}{\end{array}}
\newcommand{\bpf}{\begin{proof}}
\newcommand{\epf}{\end{proof}}

\newcommand{\Z}{\mathbb{Z}}
\newcommand{\Q}{\mathbb{Q}}
\newcommand{\Zp}{\mathbb{Z}_{p}}
\newcommand{\Qp}{\mathbb{Q}_{p}}
\newcommand{\Op}{\mathcal{O}}
\newcommand{\al}{\alpha}

\newcommand{\Ga}{\Gamma}

\newcommand{\la}{\lambda}

\newcommand{\Hi}{H_{\mathrm{Iw}}}

\DeclareMathOperator{\Gal}{Gal}
\DeclareMathOperator{\Hom}{Hom} \DeclareMathOperator{\rank}{rank}
\DeclareMathOperator{\corank}{corank}
\DeclareMathOperator{\Ext}{Ext}

\newcommand{\cyc}{\mathrm{cyc}}
\newcommand{\ord}{\mathrm{ord}}
\newcommand{\sep}{\mathrm{sep}}

\newcommand{\Iw}{\mathrm{Iw}}

\newcommand{\WK}{WK^{\acute{e}t}}

\newcommand{\cG}{\mathcal{G}}

\newcommand{\ot}{\otimes}
\newcommand{\ilim}{\displaystyle \mathop{\varinjlim}\limits}
\newcommand{\plim}{\displaystyle \mathop{\varprojlim}\limits}

\newcommand{\coker}{\mathrm{coker}\,}

\newcommand{\lra}{\longrightarrow}

\newcommand{\ps}[1]{[[ #1 ]]}

  \DeclareFontFamily{U}{wncy}{}
  \DeclareFontShape{U}{wncy}{m}{n}{<->wncyr10}{}
  \DeclareSymbolFont{mcy}{U}{wncy}{m}{n}
  \DeclareMathSymbol{\sha}{\mathord}{mcy}{"58}

\numberwithin{equation}{subsection}

\begin{document}
\title{On the codescent of \'etale wild kernels in $p$-adic Lie extensions}
 \author{Meng Fai Lim\footnote{School of Mathematics and Statistics $\&$ Hubei Key Laboratory of Mathematical Sciences, Central China Normal University, Wuhan, 430079, P.R.China. E-mail: \texttt{limmf@ccnu.edu.cn}} }
\date{}
\maketitle

\begin{abstract} \footnotesize
\noindent Let $F$ be a number field. We estimate the kernels and cokernels of the codescent maps of the \'etale wild kernels over various $p$-adic Lie extensions. For this, we propose a novel approach of viewing the \'etale wild kernel as an appropriate fine Selmer group in the sense of Coates-Sujatha. This viewpoint reduces the problem to a control theorem of the said fine Selmer groups, which in turn allows us to employ the strategies developed by Mazur and Greenberg. As applications of our estimates on the kernels and cokernels of the codescent maps, we establish asymptotic growth formulas for the \'etale wild kernels in the various said $p$-adic Lie extensions. We then relate these  growth formulas to Greenberg's conjecture (and its noncommutative analogue). Finally, we shall give some examples to illustrate our results.


\end{abstract}

\section{Introduction}

Throughout the paper, $p$ will denote a fixed prime. Let $F$ be a number field. In the event that $p=2$, assume further that $F$ has no real primes. The classical wild kernel $K_2^w(F)$ fits into the Moore exact sequence
\[ 0 \lra K_2^w(F) \lra K_2(F) \lra \mu(F_v) \lra \mu(F)\lra 0 \]
(cf.\ \cite[p.\ 157]{Milnor}). Its $p$-primary part can be described as the kernel of the localization map
\[ H^2(G_S(F), \Zp(2))\lra \bigoplus_{v\in S}H^2(F_v,\Zp(2)) \]
of continuous cohomology groups, where $S$ is a finite set of primes of $F$ containing those primes above $p$ and infinity. Inspired by this latter observation, Kolster \cite{Kol} and Nguyen Quang Do \cite{NQD,NQD02} defined higher analog of the $p$-primary part of the classical wild kernel which they coined as the \'etale wild kernels $\WK_{2i}(F)$ (see Subsection \ref{etale wild subsec} for the precise definition). %These higher \'etale wild kernels also fit into higher analogue of Moore's exact sequence and are related to the Quillen-Lichtenbaum conjecture for the (odd) torsion part of the even $K$-groups (for instance, see \cite{Ban, Ban13, KolTrieste}).
We should also remark that this \'etale wild kernel has been studied in a slightly different guise by Schneider \cite{Sch79}. These \'etale wild kernels can be thought as analogues of the ($p$-primary) ideal class groups (see \cite{Ban}) and they are central in many deep questions in algebraic number theory. For instances, the conjecture of Kummer–Vandiver can be reformulated in terms of the wild kernels (see \cite{BG}). More importantly, they can be related to the special values of the Dedekind zeta function which is an equivalent form of Lichtenbaum conjecture (see \cite{Ban, Ban13, KolTrieste, Lic}; also see \cite{Sch79}). We should mention that this latter conjecture is known when $F$ is either a totally real field or an abelian field (cf. \cite{BNQD, KNF}) building on the main conjecture proved by Wiles \cite{Wi}. Recently, Nickel has investigated annihilation problem of the \'etale wild kernel and showed that such problem are related to the equivariant Tamagawa number conjecture (see \cite{Nic}).

A curious observation is that despite \'etale wild kernel's intimate relation in Iwasawa theory (see \cite{Kol, NQD, NQD02, Nic}), the asymptotic behaviour of the \'etale wild kernels (in the spirit of \cite{Iw59}) over a $p$-adic Lie extension has not been studied before. In fact, to the best knowledge of the author, even the case of a $\Zp$-extension does not seem to be written down in literature. The aim of the paper is therefore to fill in this gap. Namely, we will establish asymptotic growth formulas for the \'etale wild kernels in a $p$-adic Lie extension $F_\infty$. For this, one is led to analysing the codescent map
\[ \Big(\plim_n \WK_{2i}(F_n)\Big)_{G_n} \lra \WK_{2i}(F_n), \]
where the $F_n$'s are certain appropriate subextensions of $F_\infty/F$ with $G_n=\Gal(F_\infty/F_n)$. We note that for a finite $p$-extension $M/F$, the kernel and cokernel of this codescent map has been well studied in \cite{Asen, AA, KM}. However, despite their explicit nature, it does not seem easy to obtain good enough estimates from these mentioned works on the kernels and cokernels for our problem in hand.

Therefore, we shall take a different approach which we explain briefly here. Firstly, via the Poitou-Tate duality, we identify the Pontryagin dual of the \'etale wild kernel with the kernel of the following localization maps
\[
H^1\big(G_S(F),\Qp/\Zp(-i)\big)\lra \bigoplus_{v\in S} H^1\big(F_v, \Qp/\Zp(-i)\big).
\]
(In fact, this is the form that Schneider \cite{Sch79} works with.) Following Coates-Sujatha \cite{CS05}, we shall call this the fine Selmer group of $\Qp/\Zp(-i)$ over $F$, which we denote by $R_i(F)$.
Under this identification, the Pontryagin dual of $\plim \WK_{2i}(L)$ turns out to be isomorphic to the kernel of
\[
H^1\big(G_S(F_\infty),\Qp/\Zp(-i)\big)\lra \bigoplus_{w\in S(F_\infty)} H^1\big(F_{\infty,w}, \Qp/\Zp(-i)\big),
\]
which will be denoted by $R_i(F_\infty)$. Therefore, the question of studying the kernel and cokernel of the codescent map is equivalent to understanding the kernel and cokernel of the restriction maps
 \[ R_i(F_n) \lra R_i(F_\infty)^{G_n}. \]
This viewpoint opens up a channel for us to mimic the techniques developed by Mazur \cite{Maz} and Greenberg \cite{Gr03} for the study of Selmer groups of abelian varieties (also see \cite{KunLim, LimConFine}).

However, we should emphasize that this improvisation is not a direct procedure.  Indeed, in tackling these sorts of problems, one is naturally led to the problem of estimating the growth of $H^k(G_n, \Qp/\Zp(i)(F_\infty))$ (for $k=1,2)$ as $n$ varies (see discussion in the beginning of Section \ref{codescent section}), where $\Qp/\Zp(-i)(F_\infty) = \big(\Qp/\Zp(-i)\big)^{G_S(F_\infty)}$. The work of Greenberg \cite{Gr03} utilized a Lie algebraic approach which is useful in showing the finiteness of these cohomology groups. However, for our purposes in hand, we also need to know the growth rate of these cohomology groups which does not seem to be accessible directly from the Lie algebraic approach. In this paper, we shall therefore study these cohomology groups directly. We succeed in obtaining estimates for these cohomology groups when the $p$-adic Lie extension is either a $\Zp^d$-extension (see Subsection \ref{Zpd subsec}), a multi-False-Tate extension (see Subsection \ref{multiFT subsec}), a $GL_2$-extension cut out by an elliptic curve without complex multiplication in (see Subsection \ref{GL2 subsec}), or a compositum of a $GL_2$-extension with multi-False-Tate extension (see Subsection \ref{GL2 FT subsec}). Interestingly, in every of these cases, the growth of the $p$-exponent of the group $H^k(G_n, \Qp/\Zp(i)(F_\infty))$ can be shown to be $O(n)$. Despite this uniform outcome, we do not have a unified elegant way of proving these simultaneously, and have to resort to a case-by-case analysis.

We shall say a little on this analysis and leave the details to Section \ref{codescent section}. For a $\Zp^d$-extension, it is well-known that the $H^k(G_n, \Qp/\Zp(i)(F_\infty))$ are finite (see \cite{Gr03, SerreLA}). Our new insight towards obtaining a growth estimate relies on an observation of Cuoco-Monsky \cite{CuoMo} that these cohomology groups have exponents bounded by $p^{dn+c}$ for some constant $c$ independent of $n$. It would seem that the utilization of such observation is absent in the work of Greenberg \cite{Gr03}. Unfortunately, this observation of Cuoco-Monsky does not carry over for noncommutative $p$-adic extensions. Thankfully, for these noncommutative extensions considered, we are able to call upon Tate's Lemma (see Lemma \ref{Tate lemma}) which is possible due to the Pontryagin dual of the \'etale wild kernel being the kernel of cohomology groups with coefficient in $\Qp/\Zp(-i)$. The utilization of Tate's Lemma is crucial for us in making headway towards obtaining estimates on the cohomology groups $H^k(G_n, \Qp/\Zp(i)(F_\infty))$ as mentioned in the previous paragraph.

As applications of these estimates, we obtain asymptotic growth formulas for the \'etale wild kernels in the various said $p$-adic Lie extensions.
We then relate the growth of \'etale wild kernels to Greenberg's conjecture. Note that Greenberg's conjecture was originally stated over the multiple $\Zp$-extension (see \cite{G01, NQD01}), but there are natural analogue formulations of it over certain noncommutative $p$-adic Lie extensions (see \cite{CS05, HSh}; also see Subsection \ref{GCC subsec}). We shall see that Greenberg's conjecture placed a constraint on the growth rate of the \'etale wild kernels in such a $p$-adic Lie extension. For some of these $p$-adic Lie extensions, we give an equivalent characterization of Greenberg's conjecture in term of the growth of \'etale wild kernels (see Corollaries \ref{Greenberg Zpd cor} and \ref{Greenberg falseT}).

We say something briefly on the situation $p=2$ and the number field $F$ has no
real primes. In this case, there are extra complications from technical
cohomological considerations. It would seem that the object to consider in this situation is the so-called positive \'etale wild kernel as developed in \cite{AAM}. We believe there should be some scope for further work in this aspect, although this will not be pursued here.

We end the introductional section giving an outline of the paper. In Section \ref{uniform group sec}, we collect preliminary facts on uniform $p$-groups and modules over their Iwasawa algebras. We also recall certain standard bounds on the cohomology groups. In Section \ref{arithmetic prelim sec}, we introduce the \'etale wild kernels and describe how they can be interpreted as appropriate fine Selmer groups in the sense of \cite{CS05, LimFine}. We also review the Iwasawa $\mu$-conjecture and Greenberg's conjecture which will be required for our subsequent discussion. Section \ref{codescent section} is where we analyse the codescent of \'etale wild kernels in $p$-adic Lie extensions. These will be applied in Section \ref{growth sec} to obtain growth formulas for the \'etale wild kernels. We then discuss the connections between these formulas and Greenberg's conjecture.
 Finally, in Section \ref{examples sec}, we give several examples to illustrate our results. Interestingly, in some of these examples, we can even obtain better and unconditional bounds on the growth of the \'etale wild kernel than those predicted by Greenberg's conjecture (see Proposition 6.2) by building on certain calculations of Sharifi \cite{Sh08}.

\subsection*{Acknowledgement}
It is a pleasure to thank John Coates, Debanjana Kundu, Antonio Lei, Andreas Nickel and Romyar Sharifi for their interest and
helpful comments. The author is also very grateful to the anonymous referee for having carefully read an earlier version of the manuscript as well as the many valuable comments and
for pointing out some inaccuracies in its previous draft. This work is supported by the National Natural Science Foundation of China under Grant No. 11771164 and the
Fundamental Research Funds for the Central Universities of CCNU under Grant No. CCNU20TD002.

\section{Uniform pro-$p$ groups} \label{uniform group sec}

We collect here several useful results on uniform pro-$p$ groups that will be required later in the paper.
For a finitely generated pro-$p$ group $G$, write $G^{p^n} = \langle g^{p^n}|~g\in G\rangle$, in other words, the group generated by the $p^n$-th-powers of elements in $G$.
We also write $G^{\{p^n\}} = \{g^{p^n}|~g\in G\}$, which is the set consisting of the $p^n$-th-powers of elements in $G$.
The pro-$p$ group $G$ is said to be powerful if $G/G^p$ is abelian.
The lower $p$-series of $G$ is defined by $P_{1}(G) = G$, and
\[ P_{n+1}(G) = \overline{P_{n}(G)^{p}[P_{n}(G),G]}, ~\mbox{for}~ n\geq 1. \]
For a powerful group $G$, it follows from \cite[Theorem 3.6]{DSMS} that $G^{p^n} =G^{\{p^n\}} = P_{n+1}(G)$.
Furthermore, the $p$-power map
\[ P_{n}(G)/P_{n+1}(G)\stackrel{\cdot p}{\longrightarrow}
P_{n+1}(G)/P_{n+2}(G)\] is surjective for every $n\geq 1$.
If the $p$-power maps are isomorphisms for all $n\geq 1$, we then say that the group $G$ is uniformly powerful (abbrv.\ uniform).
Then $[G:P_2(G)] = [P_n(G): P_{n+1}(G)]$ for every $n\geq 1$.
Consequently, one has $[G: P_{n+1}(G)] = p^{dn}$, where $d= \dim G$ (see \cite[Definition 4.1]{DSMS}).
A well-known result of Lazard (cf.\ \cite[Corollary 8.34]{DSMS}) asserts that a compact $p$-adic Lie group always contains an open normal uniform subgroup.
Therefore, one can always reduce consideration for a general compact $p$-adic Lie group to the case of a uniform group, which we will do throughout the paper.
In particular, for a uniform group, we have $G^{p^n} =G^{\{p^n\}}$. This latter fact and the following lemma will be utilized without further mention.

\bl
Let $G$ be a uniform group and $N$ a closed normal subgroup of $G$ such that $R:=G/N$ is uniform.
Then $N$ is also uniform.
Furthermore, writing $N_n= N^{p^n}$, $G_n=G^{p^n}$, and $R_n=R^{p^n}$, we have $N_n = G_n\cap N$ and $G_n/N_n\cong R_n$.
\el

\bpf
See \cite[Lemma 2.6]{HungLim}.
\epf

\subsection{Iwasawa invariants}
Let $G$ denote a uniform pro-$p$ group. The completed group algebra of $G$ over $\Zp$ is given by
 \[ \Zp\ps{G} = \plim_U \Zp[G/U], \]
where $U$ runs over the open normal subgroups of $G$ and the inverse
limit is taken with respect to the canonical projection maps. It is well known that $\Zp\ps{G}$ is
a Noetherian Auslander regular ring without zero divisors (cf.\ \cite[Theorem 3.26]{V02} or \cite[Theorem A.1]{LimFine}; also see \cite{Neu}).
Hence it admits a skew field $Q(G)$ which is flat
over $\Zp\ps{G}$ (see \cite[Chapters 6 and 10]{GW} or \cite[Chapter
4, \S 9 and \S 10]{Lam}). As a consequence, one can define the $\Zp\ps{G}$-rank of a finitely generated $\Zp\ps{G}$-module $M$ by setting
\[\rank_{\Zp\ps{G}}(M)  = \dim_{Q(G)} (Q(G)\ot_{\Zp\ps{G}}M).\]
The $\Zp\ps{G}$-module $M$ is said to be torsion  if $\rank_{\Zp\ps{G}} (M) = 0$.
It is a standard fact that $M$ is torsion over $\Zp\ps{G}$ if and only if $\Hom_{\Zp\ps{G}}(M,\Zp\ps{G})=0$ (for instance, see \cite[Lemma 4.2]{LimFine}). In the event that the torsion $\Zp\ps{G}$-module $M$ satisfies $\Ext^1_{\Zp\ps{G}}(M,\Zp\ps{G})=0$, we shall say that $M$ is a pseudo-null $\Zp\ps{G}$-module.

For a finitely generated $\Zp\ps{G}$-module $M$, denote by
$M[p^\infty]$ the $\Zp\ps{G}$-submodule of $M$ consisting of elements
of $M$ which are annihilated by some power of $p$. Howson \cite[Proposition
1.11]{Ho2}, and independently Venjakob \cite[Theorem 3.40]{V02}), showed that there is a
$\Zp\ps{G}$-homomorphism
\[ \varphi: M[p^\infty]\lra \bigoplus_{i=1}^s\Zp\ps{G}/p^{\alpha_i},\] whose
kernel and cokernel are pseudo-null $\Zp\ps{G}$-modules, and where
the integers $s$ and $\alpha_i$ are uniquely determined. The $\mu_G$-invariant of $M$ is then defined to be $\mu_G(M) = \displaystyle
\sum_{i=1}^s\alpha_i$.

\subsection{Some estimates on cohomology groups}

In this subsection, we record certain basic estimates on the cohomology groups.
For an abelian group $M$, denote by $M[p^j]$ the subgroup of $M$ consisting of elements of $M$ annihilated by $p^j$. In particular, we have $M[p^\infty] = \cup_{j\geq 1}M[p^j]$.
For a discrete $p$-primary abelian group or a compact pro-$p$ abelian group $M$, its Pontryagin dual is defined by $M^\vee = \Hom_{cont}(M, \Qp/\Zp)$, i.e., the set of continuous group homomorphisms from $M$ to $\Qp/\Zp$, where $\Qp/\Zp$ is given the discrete topology.
If $G$ is a profinite group and $M$ a $G$-module, we let $M^G$ denote the subgroup of $M$ consisting of elements fixed by $G$, and set $M_G$ to be the largest quotient of $M$ on which $G$ acts trivially.
In particular, if $M$ is a discrete $G$-module, we denote by $H^k\left(G, M\right)$ the $k$-th Galois cohomology group of $G$ with coefficients in $M$. For a finite $p$-group $N$, we write $\ord_p(N)$ for the $p$-exponent of $N$, i.e., $|N| = p^{\ord_p(N)}$.

\bl
\label{lemma: LM15 Lemma 3.2}
Let $G$ be a pro-$p$ group. Suppose that $M$ is a discrete $G$-module which is cofinitely generated over $\Zp$. Write $M_{div}$ for the maximal $p$-divisible subgroup of $M$.
If $h_1(G) = \dim_{\Z/p\Z}\left( H^1\left( G, \Z/p\Z\right)\right)$ is finite, then
\[
\dim_{\Z/p\Z}\left( H^1\left( G,  M\right)[p]\right) \leq h_1(G)\left( \corank_{\Zp} (M) + \ord_p(M/M_{div})\right).
\]
If $h_2(G) = \dim_{\Z/p\Z}\left( H^2\left( G,\ \Z/p\Z\right)\right)$ is finite, then
\[
\dim_{\Z/p\Z}\left( H^2\left( G,  M\right)[p]\right) \leq h_2(G)\left( \corank_{\Zp} (M) + \ord_p(M/M_{div})\right).
\]
\el

\begin{proof}
The first inequality is proven in \cite[Lemma 3.2]{LimMurty15}.
The second inequality is proven similarly.
\end{proof}

\bl \label{lemma: LM15 Lemma 3.2 finite variant}
Let $G$ be a pro-$p$ group.  Suppose that $M$ is a finite discrete $G$-module.
If $h_1(G) = \dim_{\Z/p\Z}\left( H^1\left( G, \Z/p\Z\right)\right)$ is finite, then $H^1\left(G,  M\right)$ is finite with
\[
\ord_p\big(H^1(G,  M)\big) \leq h_1(G)\ \ord_p (M).
\]
If $h_2(G) = \dim_{\Z/p\Z}\left( H^2\left( G,\Z/p\Z\right)\right)$ is finite, then $H^2\left( G, M\right)$ is finite with
\[
\ord_p\big(H^2(G, M)\big) \leq h_2(G)\ \ord_p (M).
\]
\el

\bpf
 This follows from a standard d\'evissage argument and noting that the only simple discrete $G$-module is $\Z/p\Z$ with trivial $G$-action (cf. \cite[Corollary 1.6.13]{NSW}).
\epf

\section{Arithmetic preliminaries} \label{arithmetic prelim sec}

\subsection{Tate twist}

Let $K$ be any field of characteristic $\neq p$. Denote by $\mu_{p^\infty}$ the group of all the $p$-power roots of unity contained in a fixed separable closure $K^{\sep}$ of $K$. The natural action of $\Gal(K^\sep/K)$ on $\mu_{p^\infty}$ induces a continuous character
\[\chi: \Gal(K^\sep/K) \lra \mathrm{Aut}(\mu_{p^\infty}) \cong \Zp^{\times}.\]
If $M$ is either a discrete or compact $\Gal(K^\sep/K)$-module, we shall write $M(i)$ for the $\Gal(K^\sep/K)$-module which is $M$ as an abelian group but with a $\Gal(K^\sep/K)$-action given by
\[ \sigma\cdot x = \chi(\sigma)^i\sigma x,  \]
where the action on the right is the original action of $\Gal(K^\sep/K)$ on $M$.
Plainly, we have $M(0)=M$ and $\mu_{p^{\infty}} \cong \Qp/\Zp(1)$. One also checks directly that
\[\Hom_{\Zp}(\Zp(i), \Qp/\Zp(j))= \Qp/\Zp(j-i) \quad \mbox{and} \quad \Hom_{\Zp}(\Qp/\Zp(i), \Qp/\Zp(j))= \Zp(j-i).\]

The continuous character
\[\chi: \Gal(K^\sep/K) \lra \mathrm{Aut}(\mu_{p^\infty}) \cong \Zp^{\times}.\]
naturally induces a group homomorphism
\[\kappa: \Gal(K(\mu_{p^{\infty}})/K) \lra \mathrm{Aut}(\mu_{p^\infty}) \cong \Zp^{\times}\]
which is sometimes called the ($p$-adic) cyclotomic character.
We now record the following well-known lemma (cf.\ \cite{Ta}) which will be frequently used in this article.

\bl[Tate's Lemma] \label{Tate lemma} Suppose that $H^0(\Gal(K^\sep/K) , \Qp/\Zp(i))$ is finite and that the Galois group $\Gal(K(\mu_{p^{\infty}})/K)$ is infinite. Then one has
\[H^1(\Gal(K(\mu_{p^{\infty}})/K) , \Qp/\Zp(i))=0.\]
\el

\bpf
 By base-changing, we may assume that $\Ga:= \Gal(K(\mu_{p^{\infty}})/K)\cong \Zp$. Then $\Qp/\Zp(i)$ can be viewed as a $\Zp\ps{\Ga}$-module. Plainly, it is cotorsion over $\Zp\ps{\Ga}$ and so it follows from \cite[Proposition 5.3.20]{NSW} that
 \[ \corank_{\Zp\ps{\Ga}}\Big(H^1(\Gal(K(\mu_{p^{\infty}})/K) , \Qp/\Zp(i))\Big) = \corank_{\Zp\ps{\Ga}}\Big(H^0(\Gal(K(\mu_{p^{\infty}})/K) , \Qp/\Zp(i))\Big).\]
 Since $H^0(\Gal(K(\mu_{p^{\infty}})/K) , \Qp/\Zp(i)) = H^0(\Gal(K^\sep/K) , \Qp/\Zp(i))$ is finite by the hypothesis, so is $H^1(\Gal(K(\mu_{p^{\infty}})/K) , \Qp/\Zp(i))$. But by \cite[Proposition 1.7.7]{NSW}, $H^1(\Gal(K(\mu_{p^{\infty}})/K) , \Qp/\Zp(i))$ is a quotient of $\Qp/\Zp(i)$ which is a $p$-divisible group, and whence, we must have
 \[H^1(\Gal(K(\mu_{p^{\infty}})/K) , \Qp/\Zp(i))=0\]
 which is what we want to show.
\epf

The lemma in particular applies when $K$ is either a finite extension of $\Qp$ or $\Q$, and $i \neq 0$.

\subsection{Iwasawa cohomology groups} \label{section Iw coh groups}

Let $F$ be a number field, and $S$ a finite set of primes of $F$ containing all the primes above $p$ and the infinite primes. Denote by $F_S$ the maximal extension of $F$ unramified outside $S$. For every extension $L$ of $F$ contained in $F_S$, we write $G_S(L) = \Gal(F_S/L)$. For every $k\geq 0$, the continuous cohomology groups $H^k(G_S(L),\Zp(i+1))$ is defined to be \[ \plim_m H^k(G_S(L),\Z/p^m\Z (i+1)).\]
Since $p$ is odd, we have $H^k(G_S(L),\Zp(i+1)) = 0$ for $k\geq 3$ (cf. \cite[Proposition 10.11.3]{NSW}).

 A Galois extension $F_\infty$ of $F$ is said to be a uniform $p$-adic Lie extension of $F$ if its Galois group $G=\Gal(F_\infty/F)$ is a uniform pro-$p$ group.
We shall always assume that our uniform $p$-adic Lie extension $F_\infty$ is contained in $F_S$ for some finite set $S$. The Iwasawa cohomology group $H^k_{\Iw, S}(F_\infty/F, \Zp(i+1))$ is then defined to be
\[H^k_{\Iw, S}(F_\infty/F, \Zp(i+1))=\plim_L H^k(G_S(L),\Zp(i+1)),\]
where $L$ runs through all the finite extensions of $F$ contained in $F_\infty$ and the transition maps are given by corestriction on cohomology. It can be shown that these cohomology groups are finitely generated over $\Zp\ps{G}$ (for instance, see \cite[Proposition 4.1.3]{LimSh}).

In this paper, we are mostly interested in the case when $i\geq 1$. Indeed, in this situation, it is by now well known that these cohomology groups are related to algebraic $K$-groups. More precisely, we have an isomorphism
\[ K_{2i}(\Op_{F,S})[p^\infty]\cong H^2(G_S(F),\Zp(i+1))\]
which is induced by the $p$-adic Chern class maps of Soul\'e \cite{Sou}. We should mention that this isomorphism is a consequence of Rost and Voevodsky \cite{Voe} (also see \cite[Section 3.2]{LimKgroups} and references therein for more details). Building on the well known fact that $H^2(G_S(L),\Zp(i+1))$ is finite for every finite extension $L$ of $F$ (cf.\ \cite{Bo, Sou}), one can even establish the following for the second Iwasawa cohomology groups.

\bl \label{torsion H2}
Assume that $i\geq 1$. Let $F_\infty$ be a uniform $p$-adic Lie extension of $F$ contained in $F_S$. Then the module
$H^2_{\Iw, S}(F_\infty/F,\Zp(i+1))$ is torsion over $\Zp\ps{\Gal(F_\infty/F)}$.
\el

\bpf
See \cite[Proposition 4.1.1]{LimKgroups}.
\epf

\br Although we are concerned with the situation $i\geq 1$, we should mention that
Schneider has conjectured that $H^2(G_S(F),\Qp/\Zp(i+1))$ is finite for $i< 0$ (see \cite[p. 192]{Sch79}). Granted this conjecture, the argument in \cite[Proposition 4.1.1]{LimKgroups} will carry over to obtain the conclusion of the preceding lemma for $i<0$.
\er

\subsection{\'Etale wild kernel} \label{etale wild subsec}

Retain the settings of Subsection \ref{section Iw coh groups}. For each $v\in S$, we have a continuous group homomorphism
 \[\Gal(\bar{F}_v/F_v) \lra \Gal(\bar{F}/F) \lra \Gal(F_S/F).\]
This in turn induces a natural map on cohomology
\[ H^k\big(G_S(F),M\big)\lra H^k\big(F_v,M\big) \]
for every $G_S(F)$-module $M$.

Let $i$ be a fixed positive integer. Following \cite{Asen, AA, Ban, Ban13, Kol, KM, NQD, NQD02, Nic}, the \'etale wild kernel $\WK_{2i}(F)$ is defined by
\[\ker\left(H^2(G_S(F),\Zp(i+1))\lra \bigoplus_{v|p} H^2(F_v, \Zp(i+1))\right).\]

By Poitou-Tate duality, the Pontryagin dual of the \'etale wild kernel fits into the following exact sequence
\[ 0 \lra \WK_{2i}(F)^{\vee} \lra H^1(G_S(F),\Qp/\Zp(-i))\lra \bigoplus_{v\in S} H^1(F_v, \Qp/\Zp(-i)). \]
In other words, the Pontryagin dual of the \'etale wild kernel can be thought as the fine Selmer group with coefficient in $\Qp/\Zp(-i)$ (in the sense of \cite{CS05, LimFine}).
We shall write $R_i(F) = \WK_{2i}(F)^{\vee}$. Note that the \'etale wild kernel is independent of the choice of $S$ (cf. \cite[Corollary 6.6]{Mi}).

Let $F_\infty$ be a uniform $p$-adic Lie extension of $F$ contained in $F_S$. We can define $R_i(L)$ and $\WK_{2i}(L)$ for every finite extension $L$ of $F$ similarly. We then write $R_i(F_\infty) = \ilim_L R_{i}(L)$ and $Y_i(F_\infty) = \plim_L\WK_{2i}(L)$, where $L$ runs through all the finite extensions of $F$ contained in $F_\infty$. Note that $R_i(F_\infty) =Y_i(F_\infty)^{\vee}$.

\bl \label{torsion Y}
Let $F_\infty$ be a uniform $p$-adic Lie extension of $F$ contained in $F_S$. Then $Y_i(F_\infty)$ is torsion over $\Zp\ps{\Gal(F_\infty/F)}$.
\el

\bpf
 Since $Y_i(F_\infty)$ is contained in $H^2_{\Iw,S}(F_\infty/F, \Zp(i+1))$, the conclusion follows from Lemma \ref{torsion H2}.
 \epf

\subsection{Iwasawa $\mu$-conjecture}

For a given $p$-adic Lie extension $\mathcal{L}$ of $F$, denote by $K(\mathcal{L})$ the maximal unramified abelian pro-$p$ extension of $\mathcal{L}$ in which every
prime above $p$ splits completely.

The cyclotomic $\Zp$-extension of $F$ will always be denoted by $F^{\cyc}$.  A well-known theorem of Iwasawa \cite{Iw73} asserts that $\Gal(K(F^\cyc)/F^\cyc)$ is a torsion $\Zp\ps{\Ga}$-module, where $\Ga=\Gal(F^\cyc/F)$. In fact, Iwasawa conjectured the following.

\begin{conjecture}[Iwasawa $\mu$-conjecture]
The group $\Gal(K(F^\cyc)/F^\cyc)$ is finitely generated over $\Zp$, or equivalently, is a torsion $\Zp\ps{\Ga}$-module with trivial $\mu_{\Ga}$-invariant.
\end{conjecture}

\br
When $F$ is an abelian extension of $\Q$, the conjecture is known to be valid by the theorem of Ferrero-Washington \cite{FW}.
\er

We now record certain consequences of the Iwasawa $\mu$-conjecture in the form of the following two lemmas. These will be required for our subsequent discussion.

\bl \label{mu=0 fg}
Let $F$ be a number field. Suppose that the Iwasawa $\mu$-conjecture holds for $F(\mu_p)^\cyc$. Then $Y_i(F^\cyc)$ is finitely generated over $\Zp$.
\el

\bpf
By \cite[Lemma 3.2]{LimFine}, it suffices to show that $Y_i\big(F(\mu_p)^\cyc\big)$ is  finitely generated over $\Zp$. Thus, we may assume that $F$ contains $\mu_p$. In this situation, the $\Zp$-finite generation of $Y_i(F^\cyc)$ is equivalent to the Iwasawa $\mu$-conjecture by \cite[Theorem 3.5]{LimFine}.
\epf

\bl \label{fg H}
Let $F_\infty$ be a pro-$p$ $p$-adic Lie extension of a number field $F$ containing $F^\cyc$. Suppose that $Y_i(F^\cyc)$ is finitely generated over $\Zp$. Then $Y_i(F_\infty)$ is finitely generated over $\Zp\ps{H}$, where $H=\Gal(F_\infty/F^\cyc)$.
\el

\bpf
See \cite[Lemma 5.2]{LimFine}.
\epf

\subsection{Greenberg's conjecture}  \label{GCC subsec}

We now come to recalling the following conjecture of Greenberg \cite[Conjecture 3.5]{G01} (also see \cite[Conjecture 4.7]{NQD01}).

\begin{conjecture}
Let $F$ be a number field and $\widetilde{F}$ the compositum of all $\Zp$-extensions of $F$. Then $\Gal(K(\widetilde{F})/\widetilde{F})$ is a pseudo-null $\Zp\ps{\Gal(\widetilde{F}/F)}$-module.
\end{conjecture}

Actually, Greenberg's conjecture is concerned with the
pseudo-nullity of a slightly bigger Galois group. For a discussion of the relation between the original conjecture of Greenberg and the slightly weaker version adopted here, we refer readers to \cite[Subsection 4.2]{LN}.

The naive noncommutative analogue of Greenberg's conjecture is false in general (see \cite[Section 5]{HSh} for counterexamples). However, if $F_\infty$ ``comes from algebraic geometry" in the sense of Fontaine-Mazur \cite{FM}, it seems quite plausible that a direct analog of Greenberg's conjecture holds, namely, $\Gal(K(F_\infty)/F_\infty)$ is a pseudo-null $\Zp\ps{\Gal(F_\infty/F)}$-module (see \cite[Question 1.3]{HSh}).

We shall primarily be interested in two classes of noncommutative $p$-adic Lie extensions in this paper. The first of which is the False-Tate extension $F(\mu_{p^{\infty}}, \al^{-p^\infty})$. Some cases of pseudo-nullity have been established by Sharifi \cite{Sh08}. We remark that Sharifi also has some positive results pertaining to the original Greenberg's conjecture for $\Q(\mu_p)$ (for $p<1000$). The other class of noncommutative $p$-adic Lie extensions that will be considered in the paper is the extension obtained by adjoining all the $p$-division points of an elliptic curve $E$ without complex multiplication. In this situation, the analogue of Greenberg's conjecture turns out to be equivalent to the Conjecture B of Coates-Sujatha which is concerned with the pseudo-nullity of the fine Selmer group of the elliptic curve $E$ in question (see \cite[Section 4]{CS05}).

For our purposes, we have the following analogous observation for the \'etale wild kernel.

\bl \label{pseudo-null same}
Suppose that $F$ contains $\mu_p$ and that $F_\infty$ is a uniform $p$-adic Lie extension containing $F^\cyc$. Writing $G=\Gal(F_\infty/F)$, we then have that $\Gal(K(F_\infty)/F_\infty)$ is a pseudo-null $\Zp\ps{G}$-module if and only if $Y_i(F_\infty)$ is a pseudo-null $\Zp\ps{G}$-module. In particular, the pseudo-nullity of
$Y_i(F_\infty)$ is independent of $i\geq 1$. \el

\bpf
We first show that every non-archimedean prime of $F_\infty$ must split completely in $K(F_\infty)$. Indeed this automatically holds for the primes above $p$ from the definition of $K(F_\infty)$. For a non-archimedean prime $w$ of $F_\infty$ outside $p$, write $v$ for the prime of $F$ below $w$. Since $v$ does not divide $p$, it follows from \cite[Theorem 7.5.3]{NSW} that the maximal pro-$p$ extension of $F_v$ is of dimension 2, whose inertia subgroup is of dimension 1. Since $F^\cyc\subseteq F_\infty$, the prime $v$ is finitely decomposed and unramified in $F^\cyc/F$. In view of the fact that $K(F_\infty)/F_\infty$ is unramified, the prime $w$ must therefore split completely in $K(F_\infty)/F_\infty$. This establishes our claim.

We now make use of the above claim to establish the following exact sequence
\[ 0\lra \Gal(K(F_\infty)/F_\infty)\lra H_{\Iw,S}^2(F_\infty/F,\Zp(1))\lra \plim_{L} \bigoplus_{w_L\in S(L)} H^2(L_{w_L},\Zp(1)). \]
where $S(L)$ denotes the set of primes of $L$ above $S$. Recall that it follows from the Poitou-Tate sequence (see \cite{NSW}) that for every finite extension $L$ of $F$ contained in $F_\infty$, we have the following exact sequence
\[ 0\lra \Gal(K_S(L)/L)\lra H^2(G_S(L),\Zp(1))\lra \plim_{L} \bigoplus_{w_L\in S(L)} H^2(L_{w_L},\Zp(1)), \]
where $K_S(L)$ is the maximal unramified abelian pro-$p$ extension of $L$ in which every prime above $S(L)$ splits completely. (Note that $\Gal(K_S(L)/L)\cong \mathrm{Cl}_S(L)[p^\infty]$, where $\mathrm{Cl}_S(L)$ is the $S$-class group of $L$.)
Taking inverse limit, we obtain the exact sequence
\[ 0\lra \Gal(K_S(F_\infty)/F_\infty)\lra H_{\Iw,S}^2(F_\infty/F,\Zp(1))\lra \plim_{L} \bigoplus_{w_L\in S(L)} H^2(L_{w_L},\Zp(1)). \]
Since $K_S(F_\infty) = K(F_\infty)$ by the claim proven in the first paragraph, this establishes the exact sequence that we require.

Since $F_\infty$ contains $\mu_{p^{\infty}}$, we may apply \cite[Lemma 2.5.1(c)]{Sh22} to conclude that
\[ \Hi^2(F_\infty/F,\Zp(i+1))\cong\Hi^2(F_\infty/F,\Zp(1))\ot\Zp(i). \]
For each $v\in S$, we fix a prime $w_v$ of $F_\infty$ above $v$ and write $G_{v}$ for the decomposition group of $G$ at this said prime. Then there is an isomorphism
\[ \plim_{L} \bigoplus_{w_L\in S(L)} H^2(L_{w_L},\Zp(j)) \cong \bigoplus_{v\in S} \mathrm{Ind}_G^{G_v}\big(\Hi^2(F_{\infty, w_v}/F_v,\Zp(j))\big)\]
for every $j$ (for instances, see \cite[Lemma 5.3.2]{LimPT}), where the local Iwaswa cohomology group is defined analogously as in the global situation and $\mathrm{Ind}_G^{G_v}$ is the compact induction in the sense of \cite[P. 737]{NSW}. Now, a similar argument to that in  \cite[Lemma 2.5.1(c)]{Sh22} yields
\[ \Hi^2(F_{\infty, w_v}/F_v,\Zp(i+1))\cong\Hi^2(F_{\infty, w_v}/F_v,\Zp(1))\ot\Zp(i). \]
In view of these observations, upon applying $-\ot\Zp(i)$ to the exact sequence obtained in the preceding paragraph, we have $Y_i(F_\infty) = \Gal(K(F_\infty)/F_\infty)(i)$. The conclusions of the lemma are now immediate from this.
\epf

 %Since $F$ contains $\mu_p$, the field $F^{\cyc}$, and hence $F_\infty$, must contain $\mu_{p^\infty}$. Therefore, the group $G_S(F_\infty)$ acts trivially on $\Qp/\Zp(-i)$ which in turn implies that \[ H^1(G_S(F_\infty),\Qp/\Zp(-i)) = \Hom(G_S(F_\infty),\Qp/\Zp(-i)). \] One has a similar conclusion for the local cohomology groups. Taking into account the above remark that every finite prime of $F_\infty$ splits completely in $K(F_\infty)$, we see that \[R_i(F_\infty) = \Hom\big(\Gal(K(F_\infty)/F_\infty),\Qp/\Zp(-i)\big).\] Upon taking Pontryagin dual, we obtain

\section{Codescent in $p$-adic extensions} \label{codescent section}
We begin setting up notation which will be adhered throughout this section without further mention.
Let $F_\infty$ be a uniform $p$-adic Lie extension of $F$ contained in $F_S$ for some appropriate set $S$ of primes. Write $G=\Gal(F_\infty/F)$ and $G_n=G^{p^n}$. The fixed field of $G_n$ is in turn denoted by $F_n$. Note that $|F_n:F|= p^{dn}$, where $d$ is the dimension of $G$.
For each $n$, there is a natural map $r_n: R_i(F_n) \lra R_i(F_\infty)^{G_n}$ induced by the restriction on cohomology. (Alternatively, one can also view this as the Pontryagin dual of the corestriction maps on cohomology.) The goal of this section is to estimate the kernel and cokernel of these maps in various $p$-adic Lie extensions. More precisely, we will carry out the analysis for a $\Zp$-extension in Subsection \ref{Zp subsec},
a $\Zp^d$-extension in Subsection \ref{Zpd subsec}, a multi-False-Tate extension in Subsection \ref{multiFT subsec}, a certain $GL_2$-extension in Subsection \ref{GL2 subsec}, and a compositum of a $GL_2$-extension with multi-False-Tate extension in Subsection \ref{GL2 FT subsec}.

We now give an overview of the general approach, leaving the details to the respective subsections. Throughout the discussion, we shall denote by $S(F_n)$ the set of primes of $F_n$ above $S$.
Consider the following commutative diagram
\[  \xymatrixrowsep{0.25in}
\xymatrixcolsep{0.15in} \entrymodifiers={!! <0pt, .8ex>+} \SelectTips{eu}{}\xymatrix{
    0 \ar[r]^{} & R_i(F_n)\ar[d]^{r_n} \ar[r] &
    H^1(G_S(F_n),\Qp/\Zp(-i))
    \ar[d]^{h_n} \ar[r] & \displaystyle\bigoplus_{v_n\in S(F_n)} H^1\left(F_{v_n}, \Qp/\Zp(-i)\right) \ar[d]^{g_n} \\
    0 \ar[r]^{} & R_i(F_\infty)^{G_n} \ar[r]^{} & H^{1}(G_S(F_\infty),\Qp/\Zp(-i))^{G_n} \ar[r] & \displaystyle\left(\bigoplus_{w\in S(F_\infty)} H^1(F_{\infty,w},\Qp/\Zp(-i))\right)^{G_n}
     } \]
with exact rows. The snake lemma yields an exact sequence
\[ 0 \lra \ker r_n \lra \ker h_n \lra C_n \lra \coker r_n \lra \coker h_n,\]
where $C_n$ is a subgroup of $\ker g_n$. Therefore, in estimating $\ker r_n$ and $\coker r_n$, one is reduced to studying $\ker h_n$, $\coker h_n$ and $\ker g_n$.
From the Hochschild-Serre spectral sequence, we see that
\[\ker h_n = H^1(G_n, \Qp/\Zp(-i)(F_\infty))\quad \mbox{and} \quad\coker h_n \subseteq H^2(G_n, \Qp/\Zp(-i)(F_\infty)),\]
 where we write $\Qp/\Zp(-i)(F_\infty) = \big(\Qp/\Zp(-i)\big)^{G_S(F_\infty)}$. This therefore leads us to estimating the groups  $H^1(G_n, \Qp/\Zp(-i)(F_\infty))$ and $H^2(G_n, \Qp/\Zp(-i)(F_\infty))$.

To analyse $\ker g_n$, we first split it into
\[\bigoplus_{v\in S}\bigoplus_{v_n|v} \ker g_{v_n},\]
 where
\[g_{v_n}:H^1\left(F_{v_n}, \Qp/\Zp(-i)\right) \lra \bigoplus_{w|v_n} H^1\left(F_{\infty, w}, \Qp/\Zp(-i)\right).\]
 Therefore, in estimating $\ker g_n$, we need to bound each $\ker g_{v_n}$ and the number of primes of $F_n$ above each $v$.
By Shapiro's lemma, we see that
\[\ker g_{v_n} = H^1\big(\Gal(F_{\infty,w}/F_{v_n}), \Qp/\Zp(-i)(F_{\infty,w})\big)\]
 for a fixed prime $w$ of $F_\infty$ above $v_n$, where $\Qp/\Zp(-i)(F_{\infty,w}) = \big(\Qp/\Zp(-i)\big)^{\Gal(\overline{F}_{\infty,w}/F_{\infty,w})}$. Thus, the problem of bounding $\ker g_{v_n}$ is the same as bounding \[H^1\big(\Gal(F_{\infty,w}/F_{v_n}), \Qp/\Zp(-i)(F_{\infty,w})\big). \]

The estimates for the cohomology groups will be dealt with in each of the respective subsections. Here, we shall take care of the estimate on the number of primes above $v\in S$. If the prime $v$ splits completely in $F_\infty/F$, then $\Gal(F_{\infty,w}/F_{v_n})=0$ and so $H^1\big(\Gal(F_{\infty,w}/F_{v_n}), \Qp/\Zp(-i)(F_{\infty,w})\big) =0$ for all $n$. It therefore suffices to consider the primes $v\in S$ which do not split completely in $F_\infty/F$. This is precisely the content of the next lemma.

 \bl \label{bound primes}
 Let $F_\infty$ be a uniform $p$-adic Lie extension of $F$ of dimension $d$. Suppose that $v$ is a prime of $F$ such that the decomposition group of $G=\Gal(F_\infty/F)$ at $v$ has dimension $t$. Then the number of primes of $F_n$ above $v$ is $O(p^{(d-t)n})$.
 \el

 \bpf
 Fix a prime of $F_\infty$ above $v$. Let $G_{n,v}$ be the decomposition group of $G_n$ at this said prime. In other words, we have $G_{n,v} = G_n\cap G_v$, where $G_v = G_{v,0}$. By \cite[Chap.\ 4, Exercise 14]{DSMS}, there exists a constant $C_v$ such that $|G_n: G_{n,v}| = C_vp^{tn}$ for sufficiently large $n$. Since $|G:G_n| = p^{dn}$, we therefore have
 \[ |G/G_n : G_v/G_{n,v}| = (1/C_v)p^{(d-t)n} =O(p^{(d-t)n}).\]
 But the index $|G/G_n : G_v/G_{n,v}|$ is precisely the number of primes of $F_n$ above $v$. Thus, this proves the lemma.
 \epf

\subsection{$\Zp$-extension} \label{Zp subsec}
As a start, we consider the case of a $\Zp$-extension. In compliance with tradition when working over a $\Zp$-extension, we shall write $\Ga=\Gal(F_\infty/F)\cong \Zp$ and $\Ga_n = \Ga^{p^n}$ here.  Since $\Zp$-extensions are unramified outside $p$ (cf.\ \cite[Theorem 1]{Iw73}), we may take $S$ in the definition of the \'etale wild kernel to consist precisely the set of primes above $p$ and infinite primes.

\bp \label{control thm Zp}
Let $i\geq 1$ be given.
The kernel and cokernel of the map $r_n: R_i(F_n) \lra R_i(F_\infty)^{\Ga_n}$ are finite and bounded independently of $n$.
\ep

\bpf We begin with the global cohomology groups
\[H^1(\Ga_n, \Qp/\Zp(-i)(F_\infty)) \quad \mbox{and}\quad H^2(\Ga_n, \Qp/\Zp(-i)(F_\infty)).\]
Since $\Ga_n\cong \Zp$, we have $H^2(\Ga_n, \Qp/\Zp(-i)(F_\infty)) =0$. Now, if $\Qp/\Zp(-i)(F_\infty)$ is finite. Since $\Ga_n$ is procyclic, it follows from \cite[Proposition 1.7.7]{NSW} that $H^1(\Ga_n, \Qp/\Zp(-i)(F_\infty)) \cong \Qp/\Zp(-i)(F_\infty)_{\Ga_n}$. Since the latter is a quotient of $\Qp/\Zp(-i)(F_\infty)$, we have
\[ \left| H^1(\Ga_n, \Qp/\Zp(-i)(F_\infty))\right|
\leq \left| \Qp/\Zp(-i)(F_\infty)\right|. \]
On the other hand, if $\Qp/\Zp(-i)(F_\infty)$ is infinite, then necessarily, $\Qp/\Zp(-i)(F_\infty) = \Qp/\Zp(-i)$. In this case, we have $H^1(\Ga_n, \Qp/\Zp(-i)) =0$ by Tate's Lemma. Either way, we see that $H^1(\Ga_n, \Qp/\Zp(-i)(F_\infty))$ is finite and bounded independently of $n$.

We now analyze $H^1\big(\Ga_{n,w}, \Qp/\Zp(-i)(F_{\infty,w})\big)$ for a fixed prime $w$ of $F_\infty$ above $v_n$. Let $v$ be the prime of $F$ below $v_n$. As seen in the discussion before Subsection \ref{Zp subsec}, it suffices to consider the primes $v\in S$ which are finitely decomposed in $F_\infty/F$. By Lemma \ref{bound primes}, the number of primes of $F_n$ above each such $v$ is bounded independently of $n$. Hence it remains to show that $H^1\big(\Ga_{n,w}, \Qp/\Zp(-i)(F_{\infty,w})\big)$ is finite and bounded independently of $n$ which in turn follows from a similar argument to that for the global cohomology group $H^1(\Ga_n, \Qp/\Zp(-i)(F_\infty))$. Hence this concludes the proof of the proposition.
\epf

\subsection{$\Zp^d$-extension} \label{Zpd subsec}

We now come to the case of a $\Zp^d$-extension, where $d\geq 2$. Since a $\Zp^d$-extension is unramified outside $p$ (cf. \cite[Theorem 1]{Iw73}), we may take $S$ to consist precisely the set of primes above $p$ and infinite primes as in the preceding subsection.

\bp \label{control thm Zpd}
Let $i\geq 1$ be given. Let $F_\infty$ be a $\Zp^d$-extension of $F$, where $d\geq 2$.
The kernel and cokernel of the map $r_n: R_i(F_n) \lra R_i(F_\infty)^{G_n}$ are finite with
\[\ord_p(\ker r_n) = O(n) \quad \mbox{and} \quad \ord_p(\coker r_n) = O(p^{(d-1)n}).\]
If the dimension of the decomposition group of $G$ at every prime of $F$ above $p$ is at least $2$, then one has $\ord_p(\coker r_n) = O(np^{(d-2)n})$.
\ep

\bpf
 Since $H^0\big(G_n,\Qp/\Zp(-i)(F_\infty)\big)$ is finite, it follows from \cite{SerreLA} that $H^k\big(G_n, \Qp/\Zp(-i)(F_\infty)\big)$ is also finite for $k=1,2$. By Lemma \ref{lemma: LM15 Lemma 3.2}, $\dim_{\Z/p\Z}H^k\big(G_n, \Qp/\Zp(-i)(F_\infty)\big)[p]$ is bounded independent of $n$ for $k=1,2$.
 On the other hand, by \cite[Lemma 2.1.1]{LiangL} (also see \cite[Theorem 2.8]{CuoMo}), there exists a constant $c$ independent of $n$ such that $p^{dn+c}$ annihilates $H^k\big(G_n, \Qp/\Zp(-i)(F_\infty)\big)$. Combining these observations, we have $\ord_p(\ker h_n) = O(n)$ and $\ord_p(\coker h_n) = O(n)$. Now, let $v\in S$ and consider
    \[\bigoplus_{v_n|v} \ker g_{v_n}=\bigoplus_{v_n|v} H^1\big(\Gal(F_{\infty,w}/F_{v_n}), \Qp/\Zp(-i)(F_{\infty,w})\big).\]
  As seen in the proof of Proposition \ref{control thm Zp}, we may assume that the prime $v$ does not decompose completely in $F_\infty/F$. If the decomposition group of $v$ is of dimension one, then a similar argument to that in the proof of Proposition \ref{control thm Zp} shows that
  \[\ord_p\Big(H^1\big(\Gal(F_{\infty,w}/F_{v_n}), \Qp/\Zp(-i)(F_{\infty,w})\big)\Big) = O(1).\] Since the decomposition group is of dimension one, the number of primes of $F_n$ above $v$ is $O(p^{(d-1)n})$ by Lemma \ref{bound primes} and so we have
  \[\ord_p\Big(\bigoplus_{v_n|v} \ker g_{v_n}\Big) = O(p^{(d-1)n}).\]
  If the decomposition group of $v$ has dimension $\geq 2$, we can apply a similar argument to that for the global cohomology groups $H^k\big(G_n, \Qp/\Zp(-i)(F_\infty)\big)$ to obtain $\ord_p(\ker g_{v_n}) = O(n)$. Since the number of primes of $F_n$ is at most $O(p^{(d-2)n})$ by Lemma \ref{bound primes}, we have \[\ord_p\Big(\bigoplus_{v_n|v} \ker g_{v_n}\Big) = O(np^{(d-2)n}).\]
  Combining these estimates, we obtain the conclusion of the proposition.
\epf

\subsection{Multi-False-Tate extensions} \label{multiFT subsec}

In this subsection, we shall always suppose that our number field $F$ contains a primitive $p$-th root of unity. Let $d\geq 2$. Consider $F_\infty= F\left(\mu_{p^\infty}, \sqrt[p^\infty]{\alpha_1} , \ldots, \sqrt[p^\infty]{\alpha_{d-1}}\right)$, where $\alpha_1, \ldots , \alpha_{d-1}\in F^{\times}$, whose image in $F^\times/(F^{\times})^p$ are linearly independent over $\Z/p\Z$.
Then $G=\Gal(F_\infty/F)\cong \Zp^{d-1}\rtimes \Zp$. Following Coates (see \cite[Section 8]{V03}), we call this a multi-False-Tate extension. Here, the set $S$ is taken to comprise precisely of all the primes of $F$ above $p$, the infinite primes and the primes that ramify in $F_\infty/F$.

\bp \label{control thm falseT}
Let $i\geq 1$ be given. Let $F_\infty$ be the multi-False-Tate extension of $F$ defined as above.
Then the kernel and cokernel of the map $r_n: R_i(F_n) \lra R_i(F_\infty)^{G_n}$ are finite with
\[\ord_p(\ker r_n) = O(n) \quad \mbox{and} \quad \ord_p(\coker r_n) = O(np^{(d-2)n}).\]
\ep

\bpf
  Write $H=\Gal(F_\infty/F(\mu_{p^\infty}))$ and $\Ga = \Gal(F(\mu_{p^\infty})/F)$.  As before, we begin by examining the kernel and cokernel of the maps $h_n$ which reduces us to estimating the cohomology groups $H^1(G_n,\Qp/\Zp(-i))$ and $H^2(G_n,\Qp/\Zp(-i))$. We shall show that
  \[ \ord_p\big(H^1(G_n,\Qp/\Zp(-i))\big) = O(n)  \quad\mbox{and}\quad \ord_p\big(H^2(G_n,\Qp/\Zp(-i))\big)= O(n)\]
  by induction on $d$. Suppose first that $d=2$. Then one has $H\cong\Ga\cong\Zp$. Since these groups have $p$-cohomological dimension one, it follows from the Hochschild-Serre spectral sequence
  \[ H^r\big(\Ga_n, H^s(H_n, \Qp/\Zp(-i))\big) \Longrightarrow H^{r+s}(G_n, \Qp/\Zp(-i))\]
  that we have a short exact sequence
  \[ 0\lra H^1\big(\Ga_n, \Qp/\Zp(-i)\big)\lra H^1(G_n,\Qp/\Zp(-i)) \lra H^1(H_n, \Qp/\Zp(-i))^{\Ga_n}  \lra 0 \]
  and an isomorphism
  \[ H^1\big(\Ga_n, H^1(H_n, \Qp/\Zp(-i))\big)\cong H^2(G_n,\Qp/\Zp(-i)).\]
  Tate's Lemma tells us that $ H^1\big(\Ga_n, \Qp/\Zp(-i)\big)=0$, and so it follows from the short exact sequence that $H^1(G_n,\Qp/\Zp(-i)) \cong H^1(H_n, \Qp/\Zp(-i))^{\Ga_n}$. Since $H_n$ acts trivially on $\Qp/\Zp(-i)$ and $H_n\cong\Zp(1)$ as $\Ga$-modules by Kummer theory, we have
  \[H^1\big(H_n, \Qp/\Zp(-i)\big)^{\Ga_n}\cong \Hom_{\Zp}\big(\Zp(1), \Qp/\Zp(-i)\big)^{\Ga_n} =  \Qp/\Zp(-1-i)^{\Ga_n} \]
  which can be seen to be finite with $\ord_p$-growth $O(n)$.
  On the other hand, we have
  \[H^2(G_n,\Qp/\Zp(-i)) \cong H^1\big(\Ga_n, H^1(H_n, \Qp/\Zp(-i))\big)\cong H^1\big(\Ga_n,\Qp/\Zp(-1-i)\big) =0,\]
  where the final zero follows from Tate's Lemma noting that $-1-i\ne 0$ for $i\geq 1$. This completes the proof of the required estimates for $d=2$. Now suppose that $d\geq 3$. We set
  \[N= \Gal\left(F_\infty/ F\left(\mu_{p^\infty}, \sqrt[p^\infty]{\alpha_1} , \ldots, \sqrt[p^\infty]{\alpha_{d-2}}\right)\right)\quad \mbox{and}\quad V = \Gal\left(F\left(\mu_{p^\infty}, \sqrt[p^\infty]{\alpha_1} , \ldots, \sqrt[p^\infty]{\alpha_{d-2}}\right)/F\right).\] Since $N_n\cong\Zp$,
  the Hochschild-Serre spectral sequence
  \[ H^r\big(V_n, H^s(N_n, \Qp/\Zp(-i))\big) \Longrightarrow H^{r+s}(G_n, \Qp/\Zp(-i))\]
   degenerates to yield an exact sequence
  \[ 0\lra H^1\big(V_n, \Qp/\Zp(-i)\big)\lra H^1(G_n,\Qp/\Zp(-i)) \lra H^1(N_n, \Qp/\Zp(-i))^{V_n}\]
  \[\lra H^2\big(V_n, \Qp/\Zp(-i)\big)\lra H^2(G_n,\Qp/\Zp(-i)) \lra H^1\big( V_n,H^1(N_n, \Qp/\Zp(-i))\big).\]
  By induction, both $H^1\big(V_n, \Qp/\Zp(-i)\big)$ and $H^2\big(V_n, \Qp/\Zp(-i)\big)$ are finite with
  \[ \ord_p\big(H^1(V_n,\Qp/\Zp(-i))\big) = O(n)  \quad\mbox{and}\quad \ord_p\big(H^2(V_n,\Qp/\Zp(-i))\big)= O(n).\]
  Now writing $U = \Gal\left(F\left(\mu_{p^\infty}, \sqrt[p^\infty]{\alpha_1} , \ldots, \sqrt[p^\infty]{\alpha_{d-2}}\right)/F(\mu_{p^\infty})\right)$, we see that $U$ acts trivially on $N$, and so the action of $V$ on $N$ factors through $\Ga$. Consequently, we have $N\cong \Zp(1)$ as $V$-modules which in turn implies that
  \[ H^1(N, \Qp/\Zp(-i)) \cong \Qp/\Zp(-1-i)\]
  as $V$-modules. Similarly, we have $H^1(N_n, \Qp/\Zp(-i)) \cong \Qp/\Zp(-1-i)$. This in turn implies that
  \[H^1(N_n, \Qp/\Zp(-i))^{V_n} \cong \big(\Qp/\Zp(-1-i)\big)^{V_n}\]
  which is finite with $\ord_p$-growth $O(n)$. On the other hand, one has
  \[H^1\big( V_n,H^1(N_n, \Qp/\Zp(-i))\big) \cong H^1\big( V_n,\Qp/\Zp(-1-i)\big)\]
  which again is finite with $\ord_p$-growth $O(n)$ by our induction hypothesis.

  It remains to consider the local restriction maps. For each prime $w$ of $F_\infty$ above $v$, we see that $F_{\infty,w}/F_v$ is a multi-False-Tate extension. Therefore, we may apply a similar argument as above to show that $\ker g_{v_n}$ is finite with $\ord_p(\ker g_{v_n}) = O(n)$. Since the decomposition group of $v\in S$ in $F_\infty/F$ has dimension $\geq 2$ (cf. \cite[Lemma 3.9]{HV}), it follows from Lemma \ref{bound primes} that the number of primes of $F_n$ above each $v$ is $O(p^{(d-2)n})$. Combining these observations, we obtain $\ord_p(\ker g_n) = O(np^{(d-2)n})$. The proof of the proposition is now completed.
\epf

\subsection{$GL_2$-extension} \label{GL2 subsec}

We now come to the case of a $GL_2$-extension cut out by the $p$-division points of an elliptic curve $E$ without complex multiplication. The following will be the main result of this subsection.

\bp \label{control thm GL2}
Let $i\geq 1$ be given. Let $E$ be an elliptic curve defined over a number field $F$ and suppose that $E$ has no complex multiplication. Suppose that $E(F)[p] = E[p]$ and that $F_\infty = F(E[p^\infty])$ is a uniform $p$-adic Lie extension of $F$. Then the kernel and cokernel of the map $r_n: R_i(F_n) \lra R_i(F_\infty)^{G_n}$ are finite with
\[\ord_p(\ker r_n) = O(n) \quad \mbox{and} \quad \ord_p(\coker r_n) = O(np^{2n}).\]
\ep

For the proof, we require three lemmas. The first of which is as follows.

\bl \label{cohomology bound nonCM}
For $j=1,2$, the group $H^j(G_n, \Qp/\Zp(-i))$ is finite with \[\ord_p\big(H^j(G_n, \Qp/\Zp(-i))\big) = O(n).\]
\el

\bpf
 Since $E$ is an elliptic curve without complex multiplication, the group $G$ has dimension $4$. By the discussion in \cite[pp. 302]{V02}, and enlarging $F$ if necessary, we may assume that $G = Z\times H$, where $Z\cong \Zp$ and $H=\Gal(F_\infty/F^{\cyc})$. Then $\Qp/\Zp(-i)^{Z_n}$ is finite, and it follows from Tate's Lemma that $H^1(Z_n, \Qp/\Zp(-i))=0$. Therefore, the spectral sequence
\[ H^r(H_n,  H^s(Z_n, \Qp/\Zp(-i))) \Longrightarrow H^{r+s}(G_n, \Qp/\Zp(-i))  \]
degenerates to yield
\[ H^j\left(H_n, \Qp/\Zp(-i)^{Z_n}\right) \cong H^j\left(G_n, \Qp/\Zp(-i)\right)\]
for every $j$.
Now, applying Lemma \ref{lemma: LM15 Lemma 3.2 finite variant},
we see that
\[ \ord_p \big(H^1(G_n, \Qp/\Zp(-i))\big) = \ord_p\left(H^1\big(H_n, \Qp/\Zp(-i)^{Z_n}\big)\right) \leq 3 ~\ord_p\big(\Qp/\Zp(-i)^{Z_n}\big) = O(n); \] \[
\ord_p \big(H^2(G_n, \Qp/\Zp(-i))\big) = \ord_p\left(H^2\big(H_n,  \Qp/\Zp(-i)^{Z_n}\big)\right) \leq 6 ~\ord_p\big(\Qp/\Zp(-i)^{Z_n}\big) = O(n).\]
where we make use of the facts that $h_1(H_n) =3$ and $h_2(H_n)=6$. This establishes the lemma.
\epf

The remaining two lemmas are concerned with estimating the local cohomology groups. As a start, we consider the case, where the prime is above $p$.

\bl \label{local cohomology bound nonCM}
Let $K$ be a finite extension of $\Qp$, and let $E$ be an elliptic curve defined over $K$ without complex multiplication. Suppose that $E(K)[p] = E[p]$ and $K_\infty = K(E[p^\infty])$ is a uniform $p$-adic Lie extension of $K$. Write $\mathcal{G} = \Gal(K_\infty/K)$.
Then the group $H^1(\cG_n, \Qp/\Zp(-i))$ is finite with \[\ord_p\big(H^1(\cG_n, \Qp/\Zp(-i))\big) = O(n).\]
\el

\bpf
Enlarging $K$ if necessary, we may assume that $E$ has no additive reduction over $K$. Now, suppose that $E$ has split multiplicative reduction.
By the theory of Tate curves, there exists $q\in K^\times$ satisfying $|q|<1$ such that $E(\overline{K}) \cong \overline{K}^{\times}/q^{\Z}$ as $\Gal(\overline{K}/K)$-modules.
Since $K$ is assumed to contain $K(E[p])$, it also contains $\mu_p$.
Therefore, $K(\mu_{p^\infty})$ is a $\Zp$-extension of $K$.
Write $H=\Gal(K_\infty/K(\mu_{p^\infty}))$ and $\Gamma= \Gal(K(\mu_{p^\infty})/K)$.
By the theory of Tate curves, $K_\infty$ is obtained from $K(\mu_{p^\infty})$ by adjoining all the $p$-power roots of $q$.
Therefore, $K_\infty$ is a False-Tate extension of $K$. The required estimate follows from a similar argument to that in the proof of Proposition \ref{control thm falseT}.

When the elliptic curve $E$ has good supersingular reduction, the dimension of $\cG$ is either 2 or 4 (cf.\ \cite[IV A.2.2]{SerreAEC}). When the dimension of $\cG$ is 4, one may now imitate the argument in the proof of Lemma \ref{cohomology bound nonCM} to prove the asserted estimate. For the dimension 2 case, we note that the Lie algebra of $\cG$ is commutative (see loc. cit.). Therefore, in this context, we may apply a similar argument to that in Proposition \ref{control thm Zpd}.

Now, suppose that $E$ has good ordinary reduction. Then the dimension of $\cG$ is 3 (cf. \cite[Proposition 2.8]{Coates99}). For our purposes, we need to delve into a more detailed analysis of this group $\cG$ following the discussion in \cite[Lemmas 2.8 and 3.15]{Coates99}. As $E$ has good ordinary reduction, there is a short exact sequence of $\cG$-modules
\[
0\longrightarrow \widehat{E}[p^\infty] \longrightarrow E[p^\infty] \longrightarrow \widetilde{E}[p^\infty]\longrightarrow 0,\]
where $\widehat{E}$ (resp., $\widetilde{E}$) denotes the formal group (resp., the maximal \'etale quotient of the $p$-divisible group) of $E$. It follows from the definition that $\Gal(\overline{K}/K^{ur})$ acts trivially on $\widetilde{E}[p^\infty]$, where $K^{ur}$ is the maximal unramified extension of $K$.
Therefore, the extension $L_\infty:=K(\widetilde{E}[p^\infty])$ is an unramified one dimensional extension of $K$ contained in $K_\infty$, which upon enlarging $K$, we may assume it to be a $\Zp$-extension.
Then $M_\infty = L_\infty(\mu_{p^\infty})$ is a $\Zp$-extension of $K(\widetilde{E}[p^\infty])$. Again, enlarging $K$, if necessary, we may assume that $\Gal(K_\infty/M_\infty)\cong \Zp$ and $\Gal(M_\infty/K)\cong \Zp^2$. Write $U= \Gal(K_\infty/M_\infty)$ and $V= \Gal(M_\infty/K)$. Plainly, $V$ acts on $U$ by conjugation, and we like to understand this action more concretely.
Let $\rho$ be the representation of $\cG$ on $T_pE$. Fix a basis of $T_pE$ which is chosen in such a way that the first element is a basis of $T_p\widehat{E}$. Under this choice of basis, the group $\cG$ can be viewed as a subgroup of $GL_2(\Zp)$ consisting of elements of the form
\[\begin{pmatrix} \eta(\sigma) &  a(\sigma)\\
                 0 &  \varepsilon(\sigma) \\ \end{pmatrix},\quad \sigma\in \cG \]
where $\eta$ (resp., $\varepsilon$) is the representation of $\cG$ on $T_p\widehat{E}$ (resp., $T_p\widetilde{E}$). Note that $a(\sigma)$ is not identically zero. Else $T_pE\ot\Qp$ will decompose into a direct sum of two one-dimensional subspaces but this contradicts the underlying assumption that $E$ has no complex multiplication (see \cite[IV A.2.4]{SerreAEC}).

Now, since the group $U$ acts trivially on both $\widetilde{E}[p^\infty]$ and $\mu_{p^{\infty}}$, it identifies with the subgroup of elements of the form
\[\begin{pmatrix} 1 &  a(\sigma)\\
                 0 &  1 \\ \end{pmatrix}, \quad \sigma\in \cG.\]
Under this identification, a straightforward calculation shows that $U\cong \Zp(\eta\varepsilon^{-1})$ as a $V$-module. Consequently, we have
\[H^1(U,\Qp/\Zp(-i)) = \Hom_{\Zp}(U, \Q/\Zp(-i)) = \Q/\Zp(\eta^{-1}\varepsilon)(-i)
=\Q/\Zp(\varepsilon^2)(-1-i)  \]
as $V$-modules, where in the final equality, we have made used of the fact that $\Qp/\Zp(\eta\varepsilon) = \Qp/\Zp(1)$ which in turn is a consequence of the Weil-pairing. One has similar conclusion for $H^1(U_n, \Qp/\Zp(-i))$. We can now estimate the cohomology group $H^1(\cG_n, \Qp/\Zp(i))$ building on the above discussion. The inflation-restriction sequence yields an exact sequence
\[0\lra H^1(V_n, \Qp/\Zp(-i)) \lra H^1(\cG_n, \Qp/\Zp(-i)) \lra  H^1(U_n, \Qp/\Zp(-i))^{V_n}. \]
Since $V\cong \Zp^2$, we can apply a similar argument to that in Proposition \ref{control thm Zpd} to conclude that $H^1(V_n, \Qp/\Zp(-i))$ is finite with $\ord_p\big(H^1(V_n, \Qp/\Zp(-i))\big) = O(n)$. On the other hand, it follows from the above discussion that
\[H^1(U_n, \Qp/\Zp(-i))^{V_n} =\big(\Q/\Zp(\varepsilon^2)(-1-i)\big)^{V_n}
\cong \big(\widetilde{E}[p^\infty]\ot\widetilde{E}[p^\infty](-1-i)\big)^{V_n},  \]
where the latter can be easily seen to be finite with $\ord_p$-growth $O(n)$. This establishes the required estimates for the case when $E$ has good ordinary reduction. The proof of the lemma is thus completed.
\epf

We now consider the primes which do not lie above $p$.

\bl \label{local cohomology bound nonCM neq p}
Let $K$ be a finite extension of $\Q_l$, where $l\neq p$. Suppose that $E$ is an elliptic curve defined over $K$ and has potential multiplicative reduction. Suppose that $E(K)[p] = E[p]$ and $K_\infty = K(E[p^\infty])$. Write $\mathcal{G} = \Gal(K_\infty/K)$.
Then $H^1(\cG_n, \Qp/\Zp(-i))$ is finite with $\ord_p\big(H^j(\cG_n, \Qp/\Zp(-i))\big) = O(n)$.
\el

\bpf
Since $E$ has potential multiplicative reduction, the group $\cG$ is of dimension 2 (cf.\ \cite[Lemma 2.8(i)]{Coates99}). By a result of Iwasawa \cite[Theorem 7.5.3]{NSW}, $K_\infty$ has no nontrivial $p$-extension, and so $H^1(K_\infty,\Qp/\Zp(-i))=0$. It now follows from this and the inflation-restriction sequence that
\[ H^1(\cG_n, \Qp/\Zp(-i)) \cong H^1(K_n, \Qp/\Zp(-i)).\]
Taking \cite[Proposition 7.3.10(ii)]{NSW} into account, we see that the latter is finite with $\ord_p$-growth $O(n)$.
\epf

We can now prove Proposition \ref{control thm GL2}.

\bpf[Proof of Proposition \ref{control thm GL2}]
Enlarging $F$ if necessary, we may assume that $E$ has no additive reduction over $F$. Therefore, we can take $S$ to be the set of primes consisting of primes above $p$ and the multiplicative primes of $E$. Lemma \ref{cohomology bound nonCM} tells us that the maps $h_n$ have finite kernel and cokernel with $\ord_p$-growth $O(n)$. We now consider the local kernels
    \[\bigoplus_{v_n|v} \ker g_{v_n}.\]
If $v$ is a prime above $p$, we see that $\ord_p(\ker g_{v_n}) = O(n)$ by Lemma \ref{local cohomology bound nonCM}. Since the decomposition group of $\cG$ at such a prime has dimension at least $2$ as seen in the proof of Lemma \ref{local cohomology bound nonCM}, we always have
\[\ord_p\left(\bigoplus_{v_n|v} \ker g_{v_n}\right) = O(np^{2n}).\]
For the primes outside $p$, since these primes are multiplicative primes by our choice of $S$, we may apply the above argument in conjunction with Lemma \ref{local cohomology bound nonCM neq p} to obtain the required local estimates. Combining these estimates, we obtain the conclusion of the proposition.
\epf

\subsection{Compositum of $GL_2$-extension and multi-False-Tate extension} \label{GL2 FT subsec}

We come to the case of the compositum of a $GL_2$-extension carved out by an elliptic curve without complex multiplication and a multi-False-Tate extension.

\bp \label{control thm GL2 falseTate}
Let $E$ be an elliptic curve defined over $F$ without complex multiplication. Suppose that $E(F)[p] = E[p]$ and that $F(E[p^\infty])$ is a uniform $p$-adic Lie extension of $F$. Let $\alpha_1, \ldots , \alpha_{r}\in F^{\times}$, whose image in $F^\times/(F^{\times})^p$ are linearly independent over $\mathbb{F}_p$. Set $F_\infty= F\left(E[p^\infty], \sqrt[p^\infty]{\alpha_1} , \ldots, \sqrt[p^\infty]{\alpha_{r}}\right)$ which is a uniform $p$-adic Lie extension of dimension $r+4$. Then the kernel and cokernel of the map $r_n: R_i(F_n) \lra R_i(F_\infty)^{G_n}$ are finite with
\[\ord_p(\ker r_n) = O(n) \quad \mbox{and} \quad \ord_p(\coker r_n) = O(np^{(r+2)n}).\]
\ep

\bpf
 We begin showing that
  \[ \ord_p\big(H^1(G_n,\Qp/\Zp(-i))\big) = O(n)  \quad\mbox{and}\quad \ord_p\big(H^2(G_n,\Qp/\Zp(-i))\big)= O(n)\]
  by induction on $r$. When $r=0$, this is Lemma \ref{cohomology bound nonCM}. Suppose that $r\geq 1$. Write \[N= \Gal\left(F_\infty/ F\left(E[p^\infty], \sqrt[p^\infty]{\alpha_1} , \ldots, \sqrt[p^\infty]{\alpha_{r-1}}\right)\right)\] and $V = \Gal\left(F\left(E[p^\infty], \sqrt[p^\infty]{\alpha_1} , \ldots, \sqrt[p^\infty]{\alpha_{r-1}}\right)/F\right)$. As $N_n\cong\Zp$, the Hochschild-Serre spectral sequence
  \[ H^r\big(V_n, H^s(N_n, \Qp/\Zp(-i))\big) \Longrightarrow H^{r+s}(G_n, \Qp/\Zp(-i))\]
 yields an exact sequence
  \[ 0\lra H^1\big(V_n, \Qp/\Zp(-i)\big)\lra H^1(G_n,\Qp/\Zp(-i)) \lra H^1(N_n, \Qp/\Zp(-i))^{V_n}\]
  \[\lra H^2\big(V_n, \Qp/\Zp(-i)\big)\lra H^2(G_n,\Qp/\Zp(-i)) \lra H^1\big( V_n,H^1(N_n, \Qp/\Zp(-i))\big).\]
  By induction, both $H^1\big(V_n, \Qp/\Zp(-i)\big)$ and $H^2\big(V_n, \Qp/\Zp(-i)\big)$ are finite with $\ord_p$-growth $O(n)$.
  Now, if we write $U = \Gal\left(F\left(E[p^\infty], \sqrt[p^\infty]{\alpha_1} , \ldots, \sqrt[p^\infty]{\alpha_{r-1}}\right)/F(\mu_{p^\infty})\right)$, we then see that $U$ acts trivially on $N$. Thus, the action of $V$ on $N$ factors through $\Ga$, and consequently, one has $N\cong \Zp(1)$ as $V$-modules. This in turn implies that
  \[ H^1(N_n, \Qp/\Zp(-i)) \cong \Qp/\Zp(-1-i)\]
  as $V$-modules. Hence we have that $H^1(N_n, \Qp/\Zp(-i))^{V_n} = \big(\Qp/\Zp(-1-i)\big)^{V_n}$ is finite with $\ord_p$-growth $O(n)$. On the other hand,
  \[H^1\big( V_n,H^1(N_n, \Qp/\Zp(-i))\big) = H^1\big( V_n,\Qp/\Zp(-1-i)\big)\]
  which again is finite with $\ord_p$-growth $O(n)$ by our induction hypothesis. Finally, the local cohomology groups can be estimated via the above induction procedure, and building on Lemmas \ref{local cohomology bound nonCM} and \ref{local cohomology bound nonCM neq p}.
\epf

\section{Growth of
\'etale wild kernels in $p$-adic Lie extensions} \label{growth sec}
In this section, we will study the growth of \'etale wild kernel in the various $p$-adic Lie extensions as considered in Section \ref{codescent section}.

\subsection{Growth in $\Zp$-extension}

\bp \label{growth thm Zp}
Let $i\geq 1$ be given. Let $F_\infty$ be a $\Zp$-extension of $F$. Then we have
\[\ord_p(\WK_i(F_n)) = \mu\big(Y_i(F_\infty)\big) p^n + \la\big(Y_i(F_\infty)\big) n +O(1).\]
\ep

\bpf
Proposition \ref{control thm Zp} tells us that $Y_i(F_\infty)_{\Ga_n}$ is finite and
\[ \left| \ord_p\big(\WK_i(F_n)\big) - \ord_p\big(Y_i(F_\infty)_{\Ga_n}\big) \right| = O(1).\]
The required conclusion follows from this and \cite[Proposition 5.3.17]{NSW}.
\epf

\subsection{Growth in $\Zp^d$-extensions}

We begin with a general result on the growth of \'etale wild kernel in a $\Zp^d$-extension.

\bp \label{growth thm Zpd}
Let $i\geq 1$ be given. Let $F_\infty$ be a $\Zp^d$-extension of $F$. Then we have
\[\ord_p\big(\WK_i(F_n)\big) = \mu_G\big(Y_i(F_\infty)\big) p^{dn} + l_0\big(Y_i(F_\infty) \big) np^{(d-1)n} +O(p^{(d-1)n}).\]
for a certain integer $l_0\big(Y_i(F_\infty)\big)\geq 0$.
\ep

\bpf
It follows from Proposition \ref{control thm Zpd} that $Y_i(F_\infty)_{G_n}$ is finite with
\[ \Big| \ord_p\big(\WK_i(F_n)\big) - \ord_p\big(Y_i(F_\infty)_{G_n}\big) \Big| = O(p^{(d-1)n}).\]
On the other hand, the finiteness of $Y_i(F_\infty)_{G_n}$ in particularly implies that $\rank_{\Zp}\big(Y_i(F_\infty)_{G_n}\big) = 0$, and therefore, the hypothesis of \cite[Theorem 3.4]{CuoMo} is satisfied. Hence we may combine the said theorem with the above estimate to obtain the required conclusion. (For the precise description of the invariant $l_0(Y_i(F_\infty))$, we refer the readers to \cite[Definition 1.2]{CuoMo}.)
\epf

We now specialize to the situation, where the $\Zp^d$-extension contains the cyclotomic $\Zp$-extension. We shall write $G=\Gal(F_\infty/F)$ and $H=\Gal(F_\infty/F^\cyc)$.

\bt \label{Greenberg Zpd}
Let $F_\infty$ be a $\Zp^d$-extension of $F$ which contains $F^{\cyc}$. Suppose that $Y_i(F^\cyc)$ is finitely generated over $\Zp$. Then $Y_i(F_\infty)$ is a pseudo-null $\Zp\ps{G}$-module if and only if
\[\ord_p\big(\WK_{2i}(F_n)\big) =  O(p^{(d-1)n}).\]
\et

\bpf
By Lemma \ref{fg H}, we have that $Y_i(F_\infty)$ is a finitely generated $\Zp\ps{H}$-module. Since each $Y_i(F_\infty)_{G_n}$ is finite, we may apply \cite[Proposition 2.4.1]{LimKgroups} to obtain
\[\ord_p\big(Y_i(F_\infty)_{G_n}\big) =   \rank_{\Zp\ps{H}}\big( Y_i(F_\infty) \big)np^{(d-1)n}+ O(p^{(d-1)n}).\]
On the other hand, again building on the fact that $Y_i(F_\infty)$ is a finitely generated $\Zp\ps{H}$-module, a well-known result of Venjakob \cite[Example 2.3 and Proposition 5.4]{V03} tells us that $Y_i(F_\infty)$ is a pseudo-null $\Zp\ps{G}$-module if and only if $Y_i(F_\infty)$ is a torsion $\Zp\ps{H}$-module. In view of the above estimate, the latter is then equivalent to saying that
\[\ord_p\big(Y_i(F_\infty)_{G_n}\big) =   O(p^{(d-1)n}).\]
Combining this with Proposition \ref{control thm Zpd}, this is the same as saying that \[\ord_p\big(\WK_{2i}(F_n)\big) =  O(p^{(d-1)n}),\]
which concludes the proof of the theorem.
\epf

\bc \label{Greenberg Zpd cor}
Assume that $\mu_p\subseteq F$. Let $\widetilde{F}$ be the compositum of all $\Zp$-extensions of $F$ and write $G=\Gal(\widetilde{F}/F)\cong\Zp^d$. Suppose that the Iwasawa $\mu$-conjecture is valid for $F^\cyc$. Then the following statements are equivalent.
\begin{enumerate}
\item[$(a)$] Greenberg's conjecture is valid. i.e., $\Gal(K(\widetilde{F})/\widetilde{F})$ is a pseudo-null $\Zp\ps{G}$-module.
\item[$(b)$] $\ord_p\big(\WK_{2i}(F_n)\big) =  O(p^{(d-1)n})$ for some $i\geq 1$.
\item[$(c)$] $\ord_p\big(\WK_{2i}(F_n)\big) =  O(p^{(d-1)n})$ for all $i\geq 1$.
\end{enumerate}
\ec

\bpf
 This follows from a combination of Lemma \ref{pseudo-null same} and Theorem \ref{Greenberg Zpd}.
\epf

\subsection{Growth in multi-False-Tate extensions}

We come to the situation of a multi-False-Tate extension  $F_\infty= F\left(\mu_{p^\infty}, \sqrt[p^\infty]{\alpha_1} , \ldots, \sqrt[p^\infty]{\alpha_{d-1}}\right)$, where $\alpha_1, \ldots , \alpha_{d-1}\in F^{\times}$ whose image in $F^\times/(F^{\times})^p$ are linearly independent over $\mathbb{F}_p$.
In this subsection, $F$ is always assumed to contain $\mu_p$ and so $\Gal(F_\infty/F)\cong \Zp^{d-1}\rtimes \Zp$. We shall write $G=\Gal(F_\infty/F)$ and $H=\Gal(F_\infty/F^\cyc)$.

\bt \label{growth thm falseT}
Let $i\geq 1$ be given. Let $F_\infty$ be the multi-False-Tate extension of $F$ as above. Suppose that the Iwasawa $\mu$-conjecture is valid for $F^\cyc/F$. Then we have
\[\ord_p(\WK_{2i}(F_n)) = \rank_{\Zp\ps{H}}\big( Y_i(F_\infty) \big)np^{(d-1)n} +O(p^{(d-1)n}).\]
\et

\bpf
 In view of the hypothesis of the theorem and Lemmas \ref{mu=0 fg} and \ref{fg H}, we have that $Y_i(F_\infty)$ is finitely generated over $\Zp\ps{H}$. Proposition \ref{control thm falseT} tells us that each $Y_i(F_\infty)_{G_n}$ is finite. It then follows from an application of \cite[Proposition 2.4.1]{LimKgroups} (of which the idea originates from \cite{Lei}) that
\[\ord_p\big(Y_i(F_\infty)_{G_n}\big) =   \rank_{\Zp\ps{H}}\big( Y_i(F_\infty) \big)np^{(d-1)n}+ O(p^{(d-1)n}).\]
The required conclusion now follows from this and Proposition \ref{control thm falseT}.
\epf

\bc \label{Greenberg falseT}
Retain the setting of Theorem \ref{growth thm falseT}.  Then the following statements are equivalent.
\begin{enumerate}
\item[$(a)$] The noncommutative analog of Greenberg's conjecture is valid for $F_\infty$. In other words, the module $\Gal(K(F_\infty)/F_\infty)$ is  pseudo-null over $\Zp\ps{G}$.
\item[$(b)$] $\ord_p\big(\WK_{2i}(F_n)\big) =  O(p^{(d-1)n})$ for some $i\geq 1$.
\item[$(c)$] $\ord_p\big(\WK_{2i}(F_n)\big) =  O(p^{(d-1)n})$ for all $i\geq 1$.
\end{enumerate}
\ec

\bpf
 This follows from a combination of Lemma \ref{pseudo-null same} and Theorem \ref{growth thm falseT}.
\epf

\subsection{Growth in $GL_2$-extension}

We now come to the $GL_2$-extension situation. Unfortunately, the structure theory for modules over such an Iwasawa algebra is not as refined as in the $\Zp^d$-case or $\Zp^r\rtimes \Zp$-case. We at least are able to rely on \cite[Proposition 2.4]{LimAsym} to obtain an asymptotic upper bound. However, since \cite[Proposition 2.4]{LimAsym} is an estimate on the quantity $\ord_p\big(M_{G_n}/p^n\big)$ for a $\Zp\ps{G}$-module $M$, we can only obtain an upper bound for $\WK_{2i}(F_n)[p^n]$ rather than the whole group $\WK_{2i}(F_n)$.

\bp \label{growth thm GL2}
Let $i\geq 1$ be given. Let $E$ be an elliptic curve defined over $F$ which has no complex multiplication. Let $F_\infty= F(E[p^\infty])$ and assume that $G=\Gal(F_\infty/F)$ is a uniform pro-$p$ group. Write $H=\Gal(F_\infty/F^\cyc)$. Suppose that the Iwasawa $\mu$-conjecture is valid for $F^\cyc$. Then we have
\[\ord_p\big(\WK_{2i}(F_n)[p^n]\big) \leq \rank_{\Zp\ps{H}}\big(Y_i(F_\infty) \big)np^{3n} +O(p^{3n}).\]
\ep

\bpf
 By Lemmas \ref{mu=0 fg} and \ref{fg H}, we have that $Y_i(F_\infty)$ is finitely generated over $\Zp\ps{H}$. Therefore, it follows from an application of \cite[Proposition 2.4]{LimAsym} that
 \[\ord_p\big(Y_i(F_\infty)_{G_n}/p^n\big) \leq \rank_{\Zp\ps{H}}\big(Y_i(F_\infty) \big)np^{3n} +O(p^{3n}).\]
 Since
 \[ \Big(Y_i(F_\infty)_{G_n}/p^n\Big)^\vee \cong \Big(Y_i(F_\infty)_{G_n}\Big)^\vee[p^n] \cong \Big(Y_i(F_\infty)^\vee\Big)^{G_n} [p^n] = R_i(F_\infty)^{G_n}[p^n], \]
 we have
 \[\ord_p\big(R_i(F_\infty)^{G_n}[p^n]\big) \leq \rank_{\Zp\ps{H}}\big(Y_i(F_\infty) \big)np^{3n} +O(p^{3n}).\]
 On the other hand, from the restriction maps of the fine Selmer groups, we have the following exact sequence
 \[ 0\lra (\ker r_n)[p^n] \lra R_i(F_n)[p^n] \lra \big(R_i(F_\infty)^{G_n}\big)[p^n].\]
 By virtue of Proposition \ref{control thm GL2}, one has $\ord_p\big((\ker r_n)[p^n]\big) = O(n)$. Combining this with the above estimate, we obtain
  $\ord_p\big(R_i(F_n)[p^n]\big) \leq \rank_{\Zp\ps{H}}\big(Y_i(F_\infty) \big)np^{3n} +O(p^{3n})$, or equivalently,
   \[ \ord_p\big(\WK_{2i}(F_n)/p^n\big) \leq \rank_{\Zp\ps{H}}\big(Y_i(F_\infty) \big)np^{3n} +O(p^{3n}). \]
   But since $\WK_{2i}(F_n)$ is finite, we have an equality $\ord_p\big(\WK_{2i}(F_n)[p^n]\big) = \ord_p\big(\WK_{2i}(F_n)/p^n\big)$ and so this proves the
proposition.
\epf

As the estimate in Proposition \ref{growth thm GL2} is not as precise as those in Theorems \ref{growth thm Zpd} and \ref{growth thm falseT}, we are not able to prove an analogue of Corollaries \ref{Greenberg Zpd cor} and \ref{Greenberg falseT}. We do have enough to at least establish the following implication.

\bc \label{Greenberg GL2}
Retain the setting of Proposition \ref{growth thm GL2}. Suppose that $\Gal(K(F_\infty)/F_\infty)$ is pseudo-null over $\Zp\ps{G}$. Then
\[\ord_p\big(\WK_{2i}(F_n)[p^n]\big) =  O(p^{3n}).\]
\ec

\bpf
By Lemma \ref{pseudo-null same}, Greenberg's conjecture is equivalent to $Y_i(F_\infty)$ being pseudo-null over $\Zp\ps{G}$. Since $Y_i(F_\infty)$ is also finitely generated over $\Zp\ps{H}$, the result of Venjakob \cite[Example 2.3 and Proposition 5.4]{V03} tells us that $Y_i(F_\infty)$ is a torsion $\Zp\ps{H}$-module, or equivalently, $\rank_{\Zp\ps{H}}\big(Y_i(F_\infty) \big)=0$. The estimate of the corollary is now immediate from this and the preceding theorem.
\epf

\subsection{Growth in compositum of $GL_2$-extension and multi-False-Tate extension}

In this subsection, we consider the case of the compositum of a $GL_2$-extension coming from an elliptic curve without complex multiplication and a multi-False-Tate extension. As the proofs are similar to those in Proposition \ref{growth thm GL2} and Corollary \ref{Greenberg GL2}, we shall merely state the conclusions and omit their proofs.

\bp \label{growth thm GL2 FalseT}
Let $E$ be an elliptic curve defined over $F$ without complex multiplication. Suppose that $E(F)[p] = E[p]$ and that $F(E[p^\infty])$ is a uniform $p$-adic Lie extension of $F$. Let $\alpha_1, \ldots , \alpha_{r}\in F^{\times}$, whose image in $F^\times/(F^{\times})^p$ are linearly independent over $\mathbb{F}_p$. Set $F_\infty= F\left(E[p^\infty], \sqrt[p^\infty]{\alpha_1} , \ldots, \sqrt[p^\infty]{\alpha_{r}}\right)$ which is a uniform $p$-adic Lie extension of dimension $r+4$. Suppose that the Iwasawa $\mu$-conjecture is valid for $F^\cyc$. Then for $i\geq 1$, we have
\[\ord_p\big(\WK_{2i}(F_n)[p^n]\big) \leq \rank_{\Zp\ps{H}}\big(Y_i(F_\infty) \big)np^{(r+3)n} +O(p^{(r+3)n}),\]
 where $H=\Gal(F_\infty/F^\cyc)$.
\ep

\bc \label{Greenberg GL2 falseT}
Retain the setting of Proposition \ref{growth thm GL2 FalseT}. Suppose that $\Gal(K(F_\infty)/F_\infty)$ is pseudo-null over $\Zp\ps{G}$. Then
\[\ord_p\big(\WK_{2i}(F_n)[p^n]\big) =  O(p^{(r+3)n}).\]
\ec

\section{Examples} \label{examples sec}
We end the paper with some examples to illustrate our results.

\begin{enumerate}
\item[$(1)$] Let $F=\Q(\mu_p)$. Now, if $p$ is an irregular prime $<1000$, a result of Sharifi \cite[Theorem 1.3]{Sh08} asserts that $\Gal(K(\widetilde{F})/\widetilde{F})$ is pseudo-null over $\Zp\ps{\Gal(\widetilde{F}/F)}$.  We can now apply Corollary \ref{Greenberg Zpd cor} to conclude that $\ord_p\big(\WK_{2i}(F_n)\big) =  O(p^{(p-3)n/2})$ for all $i\geq 1$.

\item[$(2)$] Let $F=\Q(\mu_p)$ and $F_{\infty} = \Q(\mu_{p^{\infty}}, p^{-p^{\infty}})$. If $p$ is an irregular prime $<1000$, a result of Sharifi \cite[Propositions 3.3 and 2.1a]{Sh08} tells us that $\Gal(K(F_\infty)/F_\infty)$ is pseudo-null over $\Zp\ps{\Gal(F_\infty)/F)}$. Corollary \ref{Greenberg falseT} then tells us that $\ord_p\big(\WK_{2i}(F_n)\big) =  O(p^{n})$ for all $i\geq 1$. We shall see below that we actually have a better estimate.

\item[$(3)$] Take $p=5$ and $F=\Q(\mu_5)$. Let $E$ be the elliptic curve $150A1$ of Cremona’s table which is given by
    \[ y^2 + xy = x^3-3x-3.\]
When $F_\infty$ is one of the following $5$-adic Lie extensions
\[ \Q(E[5^{\infty}], 3^{5^{-\infty}}),\ \Q(E[5^{\infty}], 2^{5^{-\infty}},3^{5^{-\infty}}),\ \Q(E[5^{\infty}], 3^{5^{-\infty}}, 5^{5^{-\infty}}),\ \Q(E[5^{\infty}], 2^{5^{-\infty}},3^{5^{-\infty}}, 5^{5^{-\infty}}), \]
the Pontryagin dual of the fine Selmer group of $E$ over such $p$-adic Lie extension is known to be pseudo-null (cf.\ \cite[Example 23]{Bh} and \cite[Section 6, Example (b)]{LimFine}) which in turn implies that $\Gal(K(F_\infty)/F_\infty)$ is pseudo-null. Hence Corollary \ref{Greenberg GL2 falseT} applies to these examples.

%\item[$(4)$] Let $p$ to be an odd prime and $F=\Q(\mu_p)$. By \cite[Example 5.3]{HSh}, there exist $\Zp$-extension $F_\infty$ of $F^\cyc$ such that $\Gal(K(F_\infty)/F_\infty)$ with arbitrary large $\Zp\ps{H}$-rank, where $H=\Gal(F_\infty/F)\cong\Zp$. From Lemma \ref{mu=0 fg} and \ref{fg H}, and the proof of Lemma \ref{pseudo-null same}, we see that the corresponding module $Y_i(F_\infty)$ has the same $\Zp\ps{H}$-rank with $\Gal(K(F_\infty)/F_\infty)$. Combining this observation with Proposition \ref{Greenberg falseT}, we see that we can always find a $\Zp$-extension $F_\infty$ of $F^\cyc$ such that the \'etale wild kernel has arbitrary large $p^n$ growth.

\end{enumerate}

We now come back to the second example and derive a better (and \textit{unconditional}) estimate than that predicted by Greenberg's conjecture. For this, we require the following preparatory lemma.

\bl \label{main alg theorem3}
 Let $G$ be a compact pro-$p$ $p$-adic Lie group which contains a closed normal subgroup $H\cong \Zp^{d-1}$ such that $G/H\cong \Zp$. Let $M$ be a $\Zp\ps{G}$-module which is finitely generated torsion over $\Zp\ps{H}$ with $\mu_H(M)=0$. Suppose further that $M_{G_n}$ is finite for every $n$. Then we have
 \[ \ord_p\big(M_{G_n}\big) = O(np^{(d-2)n}).\]
\el

\bpf
Fix a subgroup $\Ga$ of $G$ such that $\Ga$ maps isomorphically to $G/H$ under the natural quotient map $G\lra G/H$.  For every subgroup $U$ of $G$, we write $M(U)$ for the $\Zp\ps{G}$-submodule of $M$ generated by all elements of the form $(u-1)x$, where $u\in U$ and $x\in M$. From the proof of \cite[Proposition 2.4.1]{LimKgroups}, there exists an $n_0$ such that whenever $n\geq n_0$, we have
   \[ \Big|\ord_p(M_{G_n})  - (n-n_0)\rank_{\Zp}(M_{H_n})\Big| \leq  \ord_p\left(\big(M/M(\Ga_{n_0})\big)_{H_n}\right) + 2\ \ord_p\big(M_{H_{n}}[p^\infty]\big).
   \]
   By \cite[Theorem 2.4.1]{LiangL} and the hypothesis that $\mu_H(M)=0$, we have
   \[ \ord_p\big(M_{H_{n}}[p^\infty]\big) = O(np^{(d-2)n}).\]
   Since $\mu_H(M)=0$, we also have $\mu_H\big(M/M(\Ga_{n_0})\big)=0$, and so \cite[Theorem 2.4.1]{LiangL} applies to yield
   \[  \ord_p\left(\big(M/M(\Ga_{n_0})\big)_{H_n}\right) = O(np^{(d-2)n}).\]
   Finally, by a result of Harris \cite[Theorem 1.10]{Har} and the hypothesis that $M$ is torsion over $\Zp\ps{H}$, we have
   \[ (n-n_0)\rank_{\Zp}(M_{H_n}) = O(np^{d-2}).\]
   Combining all of the above estimates, we obtain the required conclusion of the lemma.
   \epf

We can now prove the following.

\bp \label{Sharifi}
 Let $F=\Q(\mu_p)$ and $F_{\infty} = \Q(\mu_{p^{\infty}}, p^{-p^{\infty}})$, where $p$ is an irregular prime $<1000$. Then
 \[\ord_p\big(\WK_{2i}(F_n)\big) =  O(n)\]
 for every $i\geq 1$.
\ep

\bpf
The results of Sharifi \cite[Propositions 3.3 and 2.1a]{Sh08} actually show that $\Gal(K(F_\infty)/F_\infty)$ is finitely generated over $\Zp$. In particular, it is torsion over $\Zp\ps{H}$ with $\mu_H(\Gal(K(F_\infty)/F_\infty))=0$, where $H= \Gal(F_\infty/F^\cyc)$. From the proof of Lemma \ref{pseudo-null same}, we see that $Y_i(F_\infty) = \Gal(K(F_\infty)/F_\infty)(i)$, and so $Y_i(F_\infty)$ is also torsion over $\Zp\ps{H}$ with trivial $\mu_H$-invariant. Proposition \ref{control thm falseT} tells us that $Y_i(F_\infty)_{G_n}$ is finite for every $n$. Hence we may apply Lemma \ref{main alg theorem3} to conclude that
\[ \ord_p\big(Y_i(F_\infty)_{G_n}\big) = O(n).\]
 Now combining this with Proposition \ref{control thm falseT} again, we have the required estimates of the proposition.
\epf

\br
The referee has raised the question of whether the stronger estimate of the preceding proposition is expected to be true in general. We shall say a bit on this aspect here. First of all, we like to mention that Sharifi's calculations rest on two components: (1) the validity of Kummer-Vandiver's conjecture and (2) certain delicate cup products computations (see \cite[Theorem 5.7, and Corollaries 5.8 and 5.9]{Sh07}). For (1), large amount of computation works have been done in verifying Kummer-Vandiver's conjecture (see \cite{HHO} and the reference therein). In view of these computations and general belief, one would expect Kummer-Vandiver's conjecture to hold (though a proof remains elusive). For the cup products computations, to the best knowledge of the author, although there has been some further computation works done, the details do not seem to have appeared in literature as yet. In conclusion, it would seem quite plausible that the estimate in Proposition \ref{Sharifi} might hold for all irregular primes $p$, although the author has to confess that he does not know of any means (other than Sharifi's approach) to tackle this problem. \er

\end{document}